\magnification=1200

\def\today{\ifcase\month\or
 January\or February\or March\or April\or May\or June\or
 July\or August\or September\or October\or November\or December\fi
 \space\number\day, \number\year}
\font\teneuf=eufm10
\font\seveneuf=eufm7
\newfam\euffam
\def\ger{\fam\euffam\teneuf}
\textfont\euffam=\teneuf
\scriptfont\euffam=\seveneuf
\scriptscriptfont\euffam=\seveneuf
\let\scr=\scriptstyle
\let\phi=\varphi
\let\epsilon=\varepsilon

\def\tf{\tilde f} \def\te{\tilde e}
\def\hb{\hfill\break}

\def\ei{\epsilon_i}
\def\eo{\epsilon_1}
\def\ej{\epsilon_j}
\def\ipp#1#2{\langle h_{#1},#2\rangle}
\def\ttf#1#2{\tf_{#1}^{a(#2)}b_{#1}}
\def\ttfo#1#2{\ttf{#1}{#2}\otimes}
\def\gg{{\ger g}}

\def\bari{\bar\imath}
\def\barj{\bar\jmath}
\def\tto#1{\buildrel#1\over\longrightarrow}
\def\otimesc{\otimes\cdots\otimes}
\def\ui{u_\infty}
\def\ol{\overline}
\def\Bi{B(\infty)}
\def\Qq{{\bf Q}(q)}
\def\uqg{U_q({\ger g})}
\def\uqgm{U_q^-({\ger g})}
\def\VL{V(\Lambda_1)}

\def\wt{\mathop{\rm wt}\nolimits}

\newdimen\depthbox
\def\textbox#1{\depthbox=\dp\strutbox\advance\depthbox by .4pt
        \ifvmode\indent\fi\setbox0\hbox{\vrule
        \vbox{\hrule\hbox{\kern1pt#1\strut\kern1pt}\hrule}%
        \vrule}\lower\depthbox\box0}
\let\tb=\textbox

\def\sbox#1{\depthbox=\dp\strutbox\advance\depthbox by .4pt    
        \setbox0\hbox{\vrule
        \vbox{\hrule
        \hbox to 10pt{\kern1pt\hfill#1\strut\kern1pt\hfill}\hrule}%
        \vrule}\lower\depthbox\box0}
\def\smbox#1{\sbox{$#1$}}
\def\sbbox#1{\depthbox=\dp\strutbox\advance\depthbox by .4pt
        \setbox0\hbox{\vrule
        \vbox {\hrule
\hbox{\kern1pt\hfill#1\hbox{\vrule height10pt depth3.5pt width 0pt}
\kern1pt\hfill}\hrule}%
        \vrule}\lower\depthbox\box0}

%
\newdimen\hoogte    \hoogte=12pt    
\newdimen\breedte   \breedte=14pt   
\newdimen\dikte     \dikte=0.5pt    
\def\beginYoung{
       \begingroup
       \def\vr{\vrule}
       \def\fbox##1{\vbox{\offinterlineskip
                    \hrule 
                    \hbox to \breedte{\strut\vr\hfill##1\hfill\vr}
                    \hrule}}
       \vbox\bgroup \offinterlineskip \tabskip=-\dikte \lineskip=-\dikte
            \halign\bgroup &\fbox{##\unskip}\unskip  \crcr }
\def\endYoung{\egroup\egroup\endgroup}
\let\by=\beginYoung \let\ey=\endYoung

\parindent=0pt
\centerline{Crystal Bases and Young Tableaux}
\footnote{}{\noindent 
 This research was supported in part by a grant from
NSERC of Canada.
\par email address: cliff@gauss.math.ualberta.ca}
\medskip\centerline{by}
\centerline{Gerald Cliff}
\bigskip
\centerline{Department of Mathematics}
\centerline{University of Alberta}
\centerline{Edmonton, Alberta}
\centerline{Canada  T6G 2G1}
\bigskip
\baselineskip=14pt
{1. \bf Introduction}
\medskip
\par 
Kashiwara [K1] has defined the notion of a crystal
basis $B(\infty)$ for the minus part $\uqgm$ of the quantized 
enveloping algebra of a semi-simple Lie algebra.
He has shown [K2] that $B(\infty)$ has a combinatorial
description given by a mapping 
$$\Psi:B(\infty)\to B(\infty)\otimes B_{i_k}\otimes B_{i_{k-1}}\otimes\cdots
  \otimes B_{i_1}$$
for a sequence $S={i_1, i_2, \ldots i_k}$ of numbers in the index set
$I$ of the simple roots, where
each $B_j$ is a certain abstract crystal, and where
the longest word $w_0$ in the Weyl group has the form
$w_0=s_{i_1}s_{i_2}\cdots s_{i_k}$.
In general the image
of the map $\Psi$ is not known.
In this paper we give an explicit description of the image of
$\Psi$, for a specific choice of $S$, in the case of classical
Lie algebras of types $A_n$, $B_n$, $C_n$, $D_n$,
using semi-standard tableaux of Kashiwara and Nakashima [KN].

\parindent=20pt

Here are our results for type $A_n$. 
We consider the map
$$\Psi:\Bi\to\Bi\otimes(B_n\otimes B_{n-1}\otimes\cdots\otimes B_1)
 \otimes(B_n\otimes B_{n-1}\otimes\cdots\otimes B_2)\otimes\cdots
 \otimes(B_n\otimes B_{n-1})\otimes B_n.$$
Elements of $B_i$ are indexed by the set $\bf Z$ of integers, and
have the form $b_i(j)$ where $j\in\bf Z$.
Let $b$ be an element of $B(\infty)$;
suppose that $\Psi(b)$ is equal to
$$u_\infty\otimes\Bigl(b_n\bigl(j_{1,n}\bigr)\otimes
b_{n-1}\bigl(j_{1,n-1}\bigr)\otimes\cdots\otimes
b_1\bigl(j_{1,1}\bigr)\Bigr)\otimes
\Bigl(b_n\bigl(j_{2,n}\bigr)\otimes
b_{n-1}\bigl(j_{2,n-1}\bigr)\otimes\cdots\otimes
b_2\bigl(j_{2,2}\bigr)\Bigr)$$
$$\otimes\cdots\otimes
\Bigl(b_n\bigl(j_{n-1,n}\bigr)\otimes
b_{n-1}\bigl(j_{n-1,n-1}\bigr)\Bigr)\otimes
b_n(j_{n,n}).\eqno(1.1)$$
Here $u_\infty$ is the unique element of $B(\infty)$ of weight 0.
For a dominant weight $\lambda$
there is a morphism $\bar\pi_\lambda$ from $\Bi$ to the crystal
basis $B(\lambda)$ of the irreducible module of highest weight
$\lambda$. Elements of $B(\lambda)$ correspond to semi-standard
Young tableaux. 
Let $T$ be the Young tableau corresponding to $\bar\pi_\lambda(b)$.
We show that for a suitably large dominant weight $\lambda$,
the integer $j_{k,l}$ appearing in (1.1) is equal to the negative
of the number of entries in the $k$-th row of $T$
which are greater than $l$. 
From this it follows that the image
of $\Psi$ consists of all elements of the form (1.1)
where the numbers $j_{k,l}$ are any integers which
satisfy
$$0\le-j_{k,n}\le -j_{k,n-1}\le \cdots \le -j_{k,k},\qquad 1\le k\le n.$$ 

For type $C_n$, we use the map $\Psi$ from $B(\infty)$ to
$B(\infty)\otimes\bigl(B_1\otimesc B_n\otimesc B_1\bigl)
\otimes\bigr(B_2\otimesc B_n\otimesc B_2\bigl)\otimesc B_n$.
We use Young tableaux of type $C$ of [KN]; these have entries which
are of the form $1,2,\ldots n,\ol n, \ol{n-1},\ldots,\ol1$. The definition
of semi-standard uses the order
$$1\prec2\prec\cdots\prec n\prec \ol n\prec\ol{n-1}\prec\cdots\prec\ol1,$$
and a condition on when $i$ and $\bari$ can occur in the same column.
Our results for type $C$ are analogous to the type $A$ results given above.
We have analogues for types $B$ and $D$; these are somewhat more
complicated than types $A$ and $C$.

\bigskip
\noindent{2. \bf Preliminaries}
\medskip
\noindent
We follow the notation in [K1], [K2], [KN]. 
For a survey of Kashiwara's theory see [K3]. Let $\ger g$ be
a finite-dimensional split semi-simple Lie algebra over $\bf Q$
having Cartan subalgebra $\ger t$, with simple
roots indexed by $I=\{1,2,\ldots, n\}$; the simple roots are
$\{\alpha_i\}_{i\in I}$, and simple coroots are
$\{h_i\}_{i\in I}$. We denote the weight lattice by
$P$, with fundamental dominant weights
$\{\Lambda_i\}_{i\in I}$. We have an inner product $(~,~)$ on
$\ger t^*$ such that $(\alpha_i,\alpha_i)$ is a non-negative integer
and set $\langle h_i,\lambda\rangle=2(\alpha_i,\lambda)/(\alpha_i,\alpha_i)$
for $i\in I$, $\lambda\in P$. The quantized enveloping algebra
$\uqg$ is the ${\bf Q}(q)$-algebra generated by $e_i$, $f_i$ ($i\in I$)
and $q^h$ ($h\in P^*)$ subject to the relations 1.1.14--1.1.18 of [K1].
We assume that $q$ is not a root of unity.
 
For a finite-dimensional $\uqg$-module $M$ and $\lambda\in P$, the
weight space $M_\lambda$ is $\{\,m\in M:q^hm=
q^{\langle h,\lambda\rangle}m$ for all $h\in P^*\}$.
We call $M$ integrable if it is the direct sum of its weight spaces.
To $M$ there is attached [K1] a {\it crystal basis} $B$, which is
a finite set of size equal to the dimension of $M$. Moreover
$B$ is the disjoint union of subsets $B_\lambda$ of size equal to
the dimension of $M_\lambda$. For $i\in I$ there are maps $\tilde e_i$
and $\tilde f_i$ from $B$ to $B\cup\{0\}$ (which reflect the decomposition of 
$M$ as a module over the copy of $U_q({\ger sl}_2)$ generated by $e_i$
and $f_i$.) For $b_1$ and $b_2$ in $B$, we have $\tf_i b_1=b_2$ if and only
if $\te_i b_2=b_1$.
For $b\in B$ we let $\epsilon_i(b)$ be the largest
integer $k$ such that $\tilde e_i^kb\ne0$ and $\phi_i(b)$ be the largest
integer $k$ such that $\tilde f_i^kb\ne 0$.
The crystal graph of $B$ is a colored oriented graph whose vertex set is
$B$, with an arrow colored by $i$ from $b_1$ to $b_2$ if $\tf_ib_1=b_2$.

If $M_1$ and $M_2$ are integrable
$\uqg$-modules with crystal bases $B_1$ and $B_2$
respectively, then $M_1\otimes M_2$ is an integrable $\uqg$-module
with crystal basis $B_1\otimes B_2$, and the action of $\tilde e_i$ and
$\tilde f_i$ on $B_1\otimes B_2$ is given as follows (cf.~ [K1], Theorem 1):

$$\eqalign{\tilde f_i(b_1\otimes b_2)&=\cases{\tilde f_ib_1\otimes b_2
 &if $\phi_i(b_1)>\epsilon_i(b_2)$,\cr
 b_1\otimes \tilde f_ib_2&if $\phi_i(b_1)\le\epsilon_i(b_2)$,\cr}\cr
\tilde e_i(b_1\otimes b_2)&=\cases{\tilde e_ib_1\otimes b_2
 &if $\phi_i(b_1)\ge\epsilon_i(b_2)$,\cr
 b_1\otimes \tilde e_ib_2
 &if $\phi_i(b_1)<\epsilon_i(b_2)$.\cr}\cr}\eqno(2.1)$$

For a dominant weight $\lambda$, $B(\lambda)$ is the crystal basis
of the irreducible $\uqg$-module $V(\lambda)$ of highest weight $\lambda$.
The unique element of $B(\lambda)_\lambda$ is denoted $u_\lambda$. Each
element of $B(\lambda)$ is of the form $\tf_{i_1}\tf_{i_2}\cdots\tf_{i_k}
u_\lambda$ for some $i_1, i_2, \ldots i_k\in I$.
The minus part $\uqgm$ is the subalgebra of $\uqg$ generated by $f_i$,
$i\in I$. It is a direct sum of its weight spaces $\uqgm_\xi=
\{\,u\in\uqgm:q^huq^{-h}=q^{\langle h,\xi\rangle}$ for all $h\in P^*\}$,
 where $\xi\in
Q_-=\{\sum n_i\alpha_i:n_i\in{\bf Z}, n_i\le0\}$.
 It has a crystal
basis $B(\infty)$ (cf.~[K1], Theorem 5) which is the disjoint union
of subsets $B(\infty)_\xi$, $\xi\in Q_-$, each of which has size equal
to the dimension of $\uqgm_\xi$. The unique element of $\Bi_0$ is denoted
$u_\infty$, and every element of $\Bi$ is of the form
 $\tf_{i_1}\tf_{i_2}\cdots\tf_{i_k}u_\infty$ 
for some $i_1, i_2, \ldots i_k\in I$.
 There are maps $\tilde f_i$ from $\Bi$ to $\Bi$ and $\te_i$
from $B(\infty)$ to $B(\infty)\cup\{0\}$.
 For a dominant weight $\lambda$ there
is a map $\bar\pi_\lambda$ from $B(\infty)$ onto $B(\lambda)\cup\{0\}$, which
maps $B(\infty)_\xi$ to $B(\lambda)_{\lambda-\xi}$, given by
$$\bar\pi_\lambda(\tf_{i_1}\tf_{i_2}\cdots\tf_{i_k}
u_\infty)=\tf_{i_1}\tf_{i_2}\cdots\tf_{i_k}
u_\lambda.$$
For each $b\in \Bi$, $\bar\pi_\lambda(b)\ne0$ for all sufficiently large
$\lambda$ (depending on $b$)
and for $\xi\in Q_-$, $\bar\pi_\lambda:B(\infty)_\xi\to
 B(\lambda)_{\lambda-\xi}$
is a bijection for sufficiently large $\lambda$.
For $b\in \Bi$, $\epsilon_i(b)=\max\{k:\tilde e_i^kb\ne0\}$, and
$\phi_i(b)$ is given by (C1) below; for sufficiently large
$\lambda$, $\ei(b)=\ei(\bar\pi_\lambda(b))$.

The crystal bases $B(\lambda)$ and $B(\infty)$ are examples of
abstract crystals [K2]. Such a crystal is a set equipped with
maps $\tilde e_i$ and $\tilde f_i$ from $B$ to $B\cup \{0\}$ where
$0$ is an ideal element not in $B$. Each $b\in B$ has a weight $\wt b$,
and there are maps $\phi_i$, $\epsilon_i$ from $B$ to 
${\bf Z}\cup\{-\infty\}$,
which satisfy axioms C1--C4 of [K2]:
\item{(C1)}$\phi_i(b)=\ei(b)+\ipp i{\wt b},\quad i\in I$;
\item{(C2)}if $b\in B$ and $\te_ib\in B$, then $\wt\te_ib=
\wt b+\alpha_i$, $\ei(\te_i b)=\ei(b)-1$, and
$\phi(\te_i b)=\phi_i(b)+1$;
\item{(C2)${}'$}if $b\in B$ and $\tf_ib\in B$, then $\wt\tf_ib=
\wt b-\alpha_i$, $\ei(\tf_i b)=\ei(b)+1$, and
$\phi(\tf_i b)=\phi_i(b)-1$;
\item{(C3)}for $b,b'\in B$, and $i\in I$, $b'=\te_ib$ if and only if 
$b=\tf_ib'$;
\item{(C4)}for $b\in B$, if $\phi_i(b)=-\infty$, then $\te_ib=\tf_ib=0$.
 
The tensor product of two crystals is
defined, and is associative; the action of $\tilde e_i$ and $\tilde f_i$
on $b_1\otimes b_2$ is given by (2.1). We have [K2, 1.3]
\vfill\eject
$$\eqalign{\epsilon_i(b_1\otimes b_2)&=\max(\epsilon_i(b_1),\epsilon_i(b_2)-
 \langle h_i,\wt b_1\rangle)\cr
          \phi_i(b_1\otimes b_2)&=\max(\phi_i(b_2),\phi_i(b_1)+
          \langle h_i,\wt b_2\rangle).\cr}\eqno(2.2)$$
Suppose that $B_1,B_2,\cdots,B_n$ are crystals,
with $b_i\in B_i$. Using (2.2) and induction, we have, essentially as in
[KN], Proposition 2.1.1,
$$\eqalign{\phi_i(b_1\otimes b_2\otimes\cdots b_n)&=
\max\bigl\{\phi_i(b_j)+\sum_{k>j}
\ipp i{\wt b_k}: j=1,2,\ldots,n\bigr\}\cr
\epsilon_i(b_1\otimes b_2\otimes\cdots b_n)&=
\max\bigl\{\epsilon_i(b_j)-\sum_{k<j}
\ipp i{\wt b_k}: j=1,2,\ldots,n\bigr\}\cr}\eqno(2.3)$$

 For $i\in I$ there is a crystal $B_i$
whose elements are $b_i(n)$, $n\in{\bf Z}$. We have
$$\displaylines{\wt b_i(n)=n\alpha_i, \qquad \tilde e_ib_i(n)=b_i(n+1),
\qquad \tilde f_ib_i(n)=b_i(n-1),\cr
\tilde e_ib_j(n)=\tilde f_ib_j(n)=0 \quad\hbox{for }i\ne j,\cr
\phi_i(b_i(n))=n,\quad\epsilon_i(b_i(n))=-n,
\quad\phi_i(b_j(n))=\epsilon_i(b_j(n))=
  -\infty\quad\hbox{for }i\ne j.\cr}$$
The element $b_i(0)$ is denoted $b_i$.

There is a notion of morphism of crystals. For each $i\in I$ there
is [K2, Theorem 2.2.1] 
a unique injective morphism from $B(\infty)$ to $B(\infty)\otimes B_i$
which commutes with each $\tilde e_j$ and $\tilde f_j$
and which sends $u_\infty$ to $u_\infty\otimes b_i$.
By iteration, for any sequence 
$S=\{i_1,i_2,\cdots i_k\}$ there is an injective
morphism from $B(\infty)$ to $B(\infty)\otimes B_{i_k}\otimesc B_{i_1}$,
and for any $b\in B(\infty)$ such a sequence $S$ can be chosen so that
$b$ is mapped to 
$u_\infty\otimes b_{i_k}\otimesc b_{i_1}$. Hence $B(\infty)$
can be considered as a subcrystal of the limit 
$B_{i_k}\otimes\cdots\otimes B_{i_1}$; in
general it is not known how to characterize this subcrystal. It is our goal
to do this, for a specific choice of sequence $S$, for $\gg$ of types
$A$, $B$, $C$, $D$.

\bigskip

\noindent{3. \bf Type A}
\medskip\noindent
In this section we assume that $\gg$ is of type $A_n$.
We recall the use of Young tableaux of type $A_n$ in [KN].
The irreducible $\uqg$-module $\VL$ (the $q$-analogue of the vector
representation of $\gg$) has dimension $n+1$ over $\Qq$; the elements
in its crystal basis $B(\Lambda_1)$ are denoted by \sbox{$i$}$\,$,
 $i\in I$.
The crystal graph of $B(\Lambda_1)$ is
$$\sbox{1}\tto1\sbox{2}\tto2\cdots\tto{n-1}\sbox{$n$}
\tto n\textbox{$n+1$}\,.$$

The irreducible $\uqg$-module $V(\Lambda_k)$
with highest weight $\Lambda_k$ for $2\le k\le n$
can be embedded in the tensor power $\VL^{\otimes k}$, and its crystal
basis $B(\Lambda_k)$ can be regarded as a subcrystal of $B(\Lambda_1)^
{\otimes k}$. By [KN, Prop.~3.3.1] $B(\Lambda_k)$ consists of
$\smbox{i_1}\otimesc\tb{$i_k$}$ with $1\le i_1<i_2<
\cdots<i_k\le n+1$.

We identify the column
$$\vcenter{\by $i_1$\cr $i_2$\cr$\vdots$\cr
 $i_k$\cr\ey} \qquad\hbox{with~~~}
\sbox{$i_1$}\otimes\sbox{$i_2$}\otimes\cdots\otimes\tb{$i_k$}\in
B(\Lambda_1)^{\otimes k}.$$

The irreducible $\uqg$-module $V(\lambda)$ with highest weight
$\lambda=\sum\lambda_i\Lambda_i$ can be embedded in
$\VL^{\otimes \lambda_1}\otimes V(\Lambda_2)^{\otimes \lambda_2}
\otimesc V(\Lambda_n)^{\otimes\lambda_n}$, so its crystal
basis $B(\lambda)$ can be viewed as a subcrystal of
$B(\Lambda_1)^{\otimes\lambda_1}\otimes B(\Lambda_2)^{\otimes\lambda_2}
\otimesc B(\Lambda_n)^{\otimes\lambda_n}$.
Thus an element $b\in B(\Lambda)$ can be identified as a tensor product
$C_1\otimesc C_p$ of columns; we write these columns right to left
(``the Japanese
writing order'' [K3, section 5])
obtaining a Young tableau.
\medskip
\noindent Example. For $\gg$ of type $A_3$, suppose $b\in B(\Lambda_1
+2\Lambda_2+\Lambda_3)$ is identified as $C_1\otimes C_2\otimes C_3\otimes C_4$
where $C_1=\smbox3\,$, $C_2=\smbox3\otimes\smbox4\,$, $C_3=\smbox2\otimes
\smbox4\,$, $C_4=\smbox1\otimes\smbox2\otimes\smbox3\,$.
Then $b$ corresponds to the Young tableau
$$\by 1 & 2 & 3 & 3\cr
      2 & 4 & 4\cr
      3\cr \ey.$$

By construction, Young tableaux have at most
$n$ rows of boxes; the length of the $i$-th row is
$\lambda_i+\lambda_{i+1}+\cdots +\lambda_n$. Such a tableau is
said to have shape $Y(\lambda)$.
The entries in the boxes are integers $m$, with $1\le m\le n+1$. 
A Young tableau is called semi-standard if the numbers in the tableau
are strictly increasing down each column and non-decreasing from
left to right along each row. The Young tableau in the example above
is semi-standard.
The Young tableau associated
to $u_\lambda$ has only $i$'s in its $i$th row.
The basic result [KN, Theorem 3.4.2] on $B(\lambda)$ is that it is the
set of all semi-standard tableaux of shape $Y(\lambda)$.

\proclaim Definition 3.1. A semi-standard Young tableau $T$ having
$r$ rows is called {\it large\/}
if  for each $i$ satisfying $1\le i<r$,
the number of $i$'s in the $i$-th row of $T$ is greater than the number
of boxes in each row below the $i$-th row.
 
In the
following example, the first tableau is large, but the second one is not.
$$\vcenter{ \by 1 & 1 & 1 &2 & 3&4\cr
               2 &2  \cr
               3\cr \ey}\qquad
 \vcenter{ \by 1 & 2 & 2 & 3&4\cr
               2 & 3 & 3 \cr
               3\cr \ey}  $$

\proclaim Lemma 3.2. Suppose that $b\in B(\infty)_\xi$, where
$\xi=-\sum n_i\alpha_i$. Suppose that $\lambda=\sum \lambda_i\Lambda_i$,
where each $\lambda_i> n_i$. Then the Young tableau given by
$\bar\pi_\lambda(b)\in B(\lambda)$ is large.

\noindent Proof: Suppose that 
$b=\tf_{i_1}\tf_{i_2}\cdots\tf_{i_k}
u_\infty$. Then $n_i$ of the indices $i_1, \ldots, i_k$ are equal to
$i$. 
Row $i$ of $u_\lambda$ has $\lambda_i+\cdots+\lambda_n$ boxes,
each of which contains $i$, and after applying $\tf_{i_1},\tf_{i_2},
\cdots\tf_{i_k}$ to $u_\lambda$,
at most $n_i$ of these $i$'s are changed to $i+1$'s. Since
$\lambda_i> n_i$, more than 
$\lambda_{i+1}+\cdots+\lambda_n$ $i$'s 
are then left in
the $i$-th row. This is greater than the number of boxes in the
next row. This completes the proof.

\proclaim Lemma 3.3. Let $b$ be the element $u_\infty\otimes \tf_n^{a(n)}b_n
\otimes \tf_{n-1}^{a(n-1)}b_{n-1}\otimesc \tf_1^{a(1)}b_1$ of the crystal
$\Bi\otimes B_n\otimes B_{n-1}\otimes\cdots\otimes B_1$,
 where $0\le a(n)\le a(n-1)\le\cdots\le a(1)$.
Then for $1\le i\le n$\hb(1) 
$\phi_i(b)=a(i-1)-a(i)$ where $a(0)=0$;
\hb(2) $\ei(b)=a(i)-a(i+1)$  where $a(n+1)=0$.

\noindent Proof: For each tensor factor $t$ of $b$,
let
$$m_i(t)=\phi_i(t)+\hbox{the sum of $\ipp i{\wt t'}$ for tensor factors
$t'$ to the right of $t$.}\eqno(3.1)$$
From (2.3), $\phi_i(b)$ is the maximum of the $m_i(t)$ as $t$ varies
over the tensor factors of $b$.
All tensor factors $t$ of $b$ have $\phi_i(t)=-\infty$ except for
$t_1=\ui$ and $t_2=\tf_i^{a(i)}b_i$.
Suppose that $1<i<n$. Then
$m_i(\tf_i^{a(i)}b_i)=-a(i)+\langle h_i,\wt \tf_{i-1}^{a(i-1)}b_{i-1}\rangle$,
and since $\wt b_{i-1}=0$, then from (C2)${}'$ in section 2,
 $\wt \tf_{i-1}^{a(i-1)}b_{i-1}=-a(i-1)\alpha_{i-1}$, so 
$$\displaylines{
\langle h_i,\wt \tf_{i-1}^{a(i-1)}b_{i-1}\rangle=
-a(i-1)\ipp i{\alpha_{i-1}}=-a(i-1)(-1)=a(i-1),\cr
m_i(t_2)=m_i(\tf_i^{a(i)}b_i)=-a(i)+a(i-1).\cr}$$
Since $\phi_i(\ui)=0$ and
 $\ipp i{\wt\tf_i^{a(i)}b_i}=-a(i)\ipp i{\alpha_i}=-2a(i)$, we have
$$\eqalign{
m_i(t_1)&=m_i(\ui)=\phi_i(\ui)+
\ipp i{\wt\tf_{i+1}^{a(i+1)}b_{i+1}}+\ipp i{\wt\tf_i^{a(i)}b_i}+
\ipp i{\wt\tf_{i-1}^{a(i-1)}b_{i-1}}\cr
&=a(i+1)-2a(i)+a(i-1)=a(i+1)-a(i)+m_i(t_2).\cr}$$
Since $a(i+1)-a(i)\le0$, then $m_i(t_1)\le m_i(t_2)$, so the maximum of
the $m_i(t)$ is $m_i(t_2)=a(i-1)-a(i)$, and (1) holds in this case.
If $i=1$ or $i=n$, the proof of (1) is similar. The proof of (2)
is also similar.

\bigskip

\proclaim Lemma A. Suppose that $T$ is a semi-standard tableau whose
first row consists entirely of 1's; then $\eo(T)=0$.

\noindent Proof: We have $\eo(T)=\max\{k:\te_1^kT\ne 0\}$.
If $\te_1T\ne0$, then $\te_1T$ is gotten from $T$ by changing a 2 to a 1.
Then $\te_1T$ would have two 1's in one of its columns, and would not
be semi-standard. However the set of semi-standard columns is stable
under the action of $\te_1$ [KN, Prop.~3.3.1]. So $\te_1T=0$, and $\eo(T)=0$.

\proclaim Proposition A.
Let $T$ be a semi-standard tableau with only one
row.
 Let $a(x)$ be the number of entries in $T$ which are $> x$.
Let $b=F(T)=u_\infty\otimes \ttf nn\otimesc\ttfo22\ttf11$ be in 
$\Bi\otimes B_n\otimes\cdots\otimes B_2\otimes B_1$.
Then \hfill\break
(1) $\phi_i(F(T))=\phi_i(T)$ for $1<i\le n$;\hfill\break
(2) $\phi_1(F(T))=\phi_1(T)\,-$ the number of entries in $T$;
\hfill\break
(3) $\ei(F(T))=\epsilon_i(T),\quad 1\le i\le n$;\hfill\break
(4) $\tf_iF(T)=F(\tf_iT)$ if $\tf_iT\ne0$.

\noindent Proof:  For $1\le i\le n$, $\phi_i(T)$ is the number
of $i$'s in $T$; if $i\ne1$, this
equals $a(i-1)-a(i)$, and this is equal to
$\phi_i(F(T))$ by Lemma 3.3(1). 
This proves (1).
The number of
1's in $T$ equals the number of entries in 
$T$ minus $a(1)$, and (2) follows from Lemma 3.3(1). Since $\ei(T)$ is the
number of $i+1$'s in $T$, (3) follows from Lemma 3.3(2).
To prove (4), $\tf_iF(T)$ is obtained from $F(T)$ by increasing $a(i)$ by
1, whereas $\tf_iT$ is obtained from $T$ by changing an $i$ to an $i+1$,
and then $F(\tf_iT)$ is obtained from $F(T)$ by increasing $a(i)$ by 1.
This completes the proof.

\proclaim Theorem A. Let $B(\infty)$ be the crystal basis of $U_q^-(\ger g)$
where $\ger g$ has type $A_n$.
For each $i$, $1\le i\le n$, define
$$B(i)=B_n\otimes B_{n-1}\otimesc B_i.$$
Let $\Psi$ be the unique morphism from $B(\infty)$ to
$\Bi\otimes B(1)\otimes B(2)\otimesc B(n)$
which maps $u_\infty$ to 
$u_\infty\otimes(b_n\otimes \cdots \otimes b_1)
\otimes(b_n\otimes \cdots \otimes b_2)
\otimesc b_n$.   
For $b\in B(\infty)$, suppose that $\lambda$ is a dominant weight such that
the tableau $T$ associated to $\bar\pi_\lambda(b)$ is large.
Let $a(i,j)$ be the number of entries in the $i$-th row of $T$ which
are $>j$; for each $i$, define
$$\beta_i=\tf_n^{a(i,n)}b_n\otimes \tf_{n-1}^{a(i,n-1)}b_{n-1}
\otimesc \tf_i^{a(i,i)}b_i\in B(i),\quad
F(T)=\ui\otimes \beta_1\otimes\beta_2\otimesc\beta_n.$$
Then $\Psi(b)=F(T)$.

\noindent Proof: We use induction on $n$.
If $n=1$, then $b$ has the form
$\tf_1^k\ui$ for some $k$. Then $k$ is
the number of $2$'s in the only row of $T$. The map $\Psi$ goes from
$\Bi$ to $\Bi\otimes B_1$ and $\Psi(b)=u_\infty\otimes\tf_1^kb_1$,
so $\Psi(b)=F(T)$. 

Now suppose that $n>1$ and that the Theorem holds
in type $A_{n-1}$.
Suppose that $b=u_\infty$. Then 
$\bar\pi_\lambda(b)=u_\lambda$, which has no numbers greater
than $i$ in its $i$-th row, and the Theorem holds for $b$.
Suppose that $b$ is an element of $B(\infty)$ such that for any
dominant weight $\lambda$ for which the tableau $T$ associated to
$\bar\pi_\lambda(b)$ is large we have
$\Psi(b)=F(T)$. 
Fix an integer $j$, $1\le j\le n$ and choose 
$\lambda$ so that both $T=\bar\pi_\lambda(b)$
and $\tf_jT$ are large; we shall prove that
$\Psi(\tf_jb)=F(\tf_jT)$.

We have $F(T)=\Psi(b)$, so for any $i$, 
 $\ei\bigl(F(T)\bigr)=\ei(\Psi(b))=\ei(b)$.
Since $\lambda$ is large, then $\ei(b)=\ei(\bar\pi_\lambda(b))=\ei(T)$.
Therefore
$$\ei\bigl(F(T)\bigr)=\ei(T)\qquad \hbox{for }1\le i\le n.\eqno(3.2)$$

Split the tableau $T$ into two parts; the left
part $T_2$ consists of all columns of length at least 2, and
the right part $T_1$ is the rest. For example,
$$\hbox{if $T$}=\vcenter{ \by 1 & 1 & 1 &2 & 3&4\cr
               2 &2  \cr
               3\cr \ey}\quad
\hbox{then }T_2=\vcenter{ \by 1 & 1\cr
               2 &2  \cr
               3\cr \ey}\,,\quad
T_1=\vcenter{\by 1&2&3&4\cr\ey}\,.$$
By our right-to-left convention, $T$ represents the tensor $T_1\otimes T_2$.
We have $\Psi(b)=\ui\otimes\beta_1\otimesc\beta_n$. Let 
$c_1=\ui\otimes\beta_1$, and $c_2=
\beta_2\otimesc\beta_n$, so $\Psi(b)=c_1\otimes c_2$.
 We claim that
$$\phi_i(T_1)>\epsilon_i(T_2)\Longleftrightarrow
\phi_i(c_1)>\epsilon_i(c_2),\qquad 1\le i\le n.\eqno(3.3)$$
Suppose that $i=1$.
Since $T$ is large, then 
$T_1$ has at least one 1, so $\phi_1(T_1)>0$.
Also, all the entries of the first row of $T_2$ are
1's, and $\eo(T_2)=0$ by Lemma A, so $\phi_1(T_1)>\eo(T_2)$.
In the notation of Proposition A, letting $a(k)=a(1,k)$, we have
$c_1=F(T_1)$, so by Proposition A(2), $\phi_1(c_1)$
is a finite number, whereas $\epsilon_1(c_2)=-\infty$, since all tensor
factors of $c_2$ are different from $B_1$. It follows that
$\phi_1(c_1)>\epsilon_1(c_2)$, and (3.3) holds if $i=1$; indeed,
$$\phi_1(T_1)>\epsilon_1(T_2),
\qquad\phi_1(c_1)>\epsilon_1(c_2).\eqno(3.4)$$
Now assume that $i>1$. We claim that
$$\phi_i(T_1)=\phi_i(c_1)\quad\hbox{and}\quad
\epsilon_i(T_2)=\epsilon_i(c_2).\eqno(3.5)$$
The first equality follows from Proposition A(1).
Let $T_2'$ be the tableau obtained by deleting the first row of $T_2$.
This deleted row consists entirely of 1's, and these 1's are irrelevant
as far as $\epsilon_i(T_2)$ is concerned since $i>1$. Thus
$\epsilon_i(T_2)=\epsilon_i(T_2')$. Now $T_2'$ is a large semi-standard
tableau involving the symbols $2,3,\ldots, n$, and therefore represents
an element of the crystal basis, for type $A_{n-1}$, of $B_{A_{n-1}}(\lambda')$
for some dominant weight $\lambda'$ of type $A_{n-1}$.
We have $T_2'=\bar\pi_{\lambda'}(b')$ for some element $b'\in B_{A_{n-1}}
(\infty)$. By induction, letting $F_{n-1}$ and $\Psi_{n-1}$ be the
analogues of $F$ and $\Psi$ (respectively) in type
$A_{n-1}$, we have $F_{n-1}(T_2')=\Psi_{n-1}(b')=\ui\otimes c_2$. From (3.2),
$$\epsilon_i(T_2')=\epsilon_i(\ui\otimes c_2)=\max\bigl(\epsilon_i(\ui),
\epsilon_i(c_2)-\ipp i{\wt\ui}\bigr)=\max\bigl(0,\ei(c_2)\bigr)=\ei(c_2).$$
Thus $\ei(T_2)=\ei(c_2)$
and (3.5) holds, and so does (3.3).
 
Suppose that $\phi_j(T_1)>\epsilon_j(T_2)$. Then $\tf_jT=
\tf_jT_1\otimes T_2$, and since $\phi_j(c_1)>\ej(c_2)$
by (3.3), $\tf_j(c_1\otimes c_2)=\tf_jc_1\otimes c_2$.
By Proposition A(4), $\tf_jc_1=\tf_jF(T_1)=F(\tf_jT_1)$, so
$$\Psi(\tf_jb)=\tf_j\Psi(b)=\tf_jF(T)=\tf_jc_1\otimes c_2=
 F(\tf_jT_1)\otimes c_2=F(\tf_jT)$$
as desired.

Suppose that $\phi_j(T_1)\le\epsilon_j(T_2)$. Then $j>1$ by (3.4).
We have $\tf_jT=T_1\otimes\tf_jT_2$;
also, from (3.3), $\tf_j\Psi(b)=c_1\otimes\tf_jc_2$.
We want to show that $F(T_1\otimes\tf_jT_2)=c_1\otimes\tf_jc_2$.
The tableau $T_1\otimes\tf_jT_2$ differs from $T_1\otimes T_2$
only below the first row.
Since the Theorem holds for $A_{n-1}$ by induction, we have
$$\Psi_{n-1}(b')=F_{n-1}(T_2')\quad\hbox{and}\quad 
\Psi_{n-1}(\tf_jb')=F_{n-1}(\tf_jT_2').$$
Then $F_{n-1}(\tf_jT_2')=\Psi_{n-1}(\tf_jb')=\tf_j\Psi_{n-1}(b')=
\tf_jF_{n-1}(T_2')=
\tf_j(\ui\otimes c_2)=\ui\otimes\tf_j c_2$. Thus
$F_{n-1}(\tf_jT_2')=\ui\otimes\tf_jc_2$,
so $F(T_1\otimes\tf_jT_2)=c_1\otimes\tf_jc_2$.
Therefore $\Psi(\tf_jb)=F(\tf_jT)$. This completes the proof.

\proclaim Corollary A.
The image of $\Psi$ consists of all elements of the form 
$\ui\otimes \beta_1\otimes\beta_2\otimesc\beta_n$ 
where
$$\beta_i=\tf_n^{a(i,n)}b_n\otimes \tf_{n-1}^{a(i,n-1)}b_{n-1}
\otimesc \tf_i^{a(i,i)}b_i\in B(i)$$
and $\{a(i,j)\}$ are any integers such that 
$$0\le a(i,n)\le a(i,n-1)\le\cdots\le a(i,i), \quad 1\le i\le n.\eqno(*)$$

\noindent Proof: For $b\in B(\infty)$, $\Psi(b)=F(T)$
where $T$ is a large semi-standard tableau, and since
$a(i,j)$ is the number of entries of the $i$-th row of $T$ which
are $>j$, $(*)$ holds. On the other hand, if $\{a(i,j)\}$ is a set
of integers for which $(*)$ holds, form a large semi-standard tableau
$T$ by having $a(i,j)$ be the number of entries in the $i$-th row
which are $>j$. Then $T$ is in the crystal basis $B(\lambda)$ for
some dominant weight $\lambda$, so $T=\bar\pi_\lambda(b)$ for
some $b\in B(\infty)$, and then $\Psi(b)=F(T)$. The Corollary is
proved.

\bigskip

\noindent{4. \bf Type C}

In this section we assume that $\gg$ is of type $C_n$. Again we follow
the notation of [KN]. We let $\alpha_n$ be the long root, and $\alpha_i$,
$1\le i\le n-1$ be the short roots. We have $\langle h_i, \alpha_{i+1}
\rangle=-1$ for $1\le i\le n-2$, $\langle h_{n-1},\alpha_n\rangle=-2$,
$\langle h_i,\alpha_j\rangle=0$ if $|i-j|>1$.
The irreducible representation $V(\Lambda_1)$
of highest weight $\Lambda_1$ has dimension $2n$ over
$\Qq$, and the elements of its crystal basis are \sbox{$i$}$\,$ and 
\sbox{$\ol \imath$}$\,$, $1\le i\le n$. 
 The crystal graph of $B(\Lambda_1)$ is
$$\sbox{1}\tto1\sbox{2}\tto2\cdots\tto{n-1}\sbox{$n$}
\tto n\sbox{$\ol n$}\tto{n-1}\cdots\tto2
\sbox{\lower1.5pt\hbox{$\ol2$}}\tto1
\sbox{\lower1pt\hbox{$\ol1$}}\,.$$
We use the order $\prec$ where 
$$1\prec2\prec\cdots\prec n\prec \ol n\prec\ol{n-1}\prec\cdots\prec\ol1.$$

Analagous to the type $A$ case, the irreducible $\uqg$-module $V(\Lambda_k)$
with highest weight $\Lambda_k$ for $2\le k\le n$
can be embedded in the tensor power $\VL^{\otimes k}$, and its crystal
basis $B(\Lambda_k)$ can be regarded as a subcrystal of $B(\Lambda_1)^
{\otimes k}$. We again identify the column
$$\vcenter{\by $i_1$\cr $i_2$\cr$\vdots$\cr
 $i_k$\cr\ey} \qquad\hbox{with~~~}
\sbox{$i_1$}\otimes\sbox{$i_2$}\otimes\cdots\otimes\tb{$i_k$}\in
B(\Lambda_1)^{\otimes k},$$
where each $i_j$ is one of $m$ or $\ol m$, $1\le m\le n$.
Define $I_k^{(C)}$ to be the set of
all columns as above, where $1\preceq i_1\prec\cdots\prec i_k\preceq\ol1$,
and where the following condition holds:
$$\hbox{ if $i_j=p$ and $i_l=\ol p$, then}\qquad 
j+(k-l+1)\le p.\eqno(4.1)$$ 
Note that this condition implies that $1$ and $\bar 1$ cannot both
be in the same column, since $j\ge1$ and $k-l+1\ge1$.
From [KN, Prop.~4.3.2] the crystal basis $B(\Lambda_k)$
is equal to $I_k^{(C)}$.

As in type $A$, the irreducible $\uqg$-module $V(\lambda)$ with highest weight
$\lambda=\sum\lambda_i\Lambda_i$ can be embedded in
$\VL^{\otimes \lambda_1}\otimes V(\Lambda_2)^{\otimes \lambda_2}
\otimesc V(\Lambda_n)^{\otimes\lambda_n}$, so its crystal
basis $B(\lambda)$ can be viewed as a subcrystal of
$B(\Lambda_1)^{\otimes\lambda_1}\otimes B(\Lambda_2)^{\otimes\lambda_2}
\otimesc B(\Lambda_n)^{\otimes\lambda_n}$.
Then $b\in B(\lambda)$ is identified with a Young tableau, using
the right to left convention on columns as in type $A$; the columns are in 
$I_k^{(C)}$. 

We will call a tableau $T$ 
{\it almost semi-standard\/} if its columns
are in $I_k^{(C)}$ ($1\le k\le n$)
and the entries in the rows are non-decreasing, from left to right,
in the order $\prec$. Such a tableau is called large as in Definition 3.1. 

\proclaim Lemma 4.1. Suppose that $T$ is a large
almost semi-standard tableau. Then
$p$ and $\ol p$ cannot occur in the same column, 
for any $p$, $1\le p\le n$.

\noindent Proof: Suppose that $p$ and $\ol p$ occur in the same
column, with $p$ in row $j$ and $\ol p$ in row $l$. Suppose that this
column has length $k$.
Since $T$ is large, the number of $j$'s in row $j$ is greater than
the number of entries in row $l$, so $p=j$.
The condition $j+(k-l+1)\le p$ of (4.1)
fails, since $j=p$ and $k-l+1\ge 1$.
So $p$ and $\ol p$ cannot be in the same column of a large semi-standard
tableau.
\medskip
A semi-standard $C$-tableau is by definition [KN, p.~317] one
which is almost semi-standard, and for which a certain condition
[KN, M.N.2] holds concerning successive columns where
$i$ and $\ol\imath$ occur in one of the columns.
By Lemma 4.1, $i$ and $\ol\imath$ cannot occur in the same column of
a large semi-standard tableau, and  
we may dispense with [KN, M.N.2] for a large
semi-standard tableau. 
The basic result for type $C$ [KN, Theorem 4.5.1] is that $B(\lambda)$
is the set of semi-standard $C$-tableaux of shape $V(\lambda)$.
Lemma 3.2 still holds, with the same proof.

\proclaim Lemma 4.2. Let $b$ be the element 
$u_\infty\otimes \tf_1^{a(\ol2)}b_1\otimes \tf_2^{a(\ol3)}b_2
\otimes\cdots \otimes\tf_{n-1}^{a(\ol n)}b_{n-1}
\otimes\tf_n^{a(n)}b_n\otimes\tf_{n-1}^{a({n-1})}b_{n-1}
\otimesc\tf_2^{a(2)}b_2
\otimes\tf_1^{a(1)}b_1$ of the crystal
$\Bi\otimes B_1\otimes\cdots\otimes B_{n-1}\otimes B_n\otimes B_{n-1}
\otimesc B_1$,
 where $0\le a(\ol2)\le\cdots\le a(\ol n)
\le a(n)\le a({n-1})\le\cdots\le a(1)$, $n\ge2$. Let $a(0)=0$.
Then
\hb(1) $\phi_i(b)=\max\bigl(-a(\ol{i+1})+a(\ol{i+2})+a(i+1)-2a(i)+a(i-1),
a(i-1)-a(i)\bigr)$, $1\le i<n-1$;
\hb(2) $\phi_{n-1}(b)=\max\bigl(-a(\ol n)+2a(n)-2a(n-1)+a(n-2),
a(n-2)-a(n-1)\bigr)$;
\hb(3) $\phi_n(b)=a(n-1)-a(n)$. 

\noindent Proof: As in the proof of Lemma 3.3,
find the maximum of $m_i(t)$, defined in equation (3.1), 
 over the tensor factors $t$ of $b$.
First, assume that $1<i<n-1$.
The only tensor factors
$t$ of $b$
for which $\phi_i(t)\ne-\infty$, reading from right
to left, are
$t_1= \tf_i^{a(i)}b_i$, 
$t_2=\tf_i^{a(\ol{i+1})}b_i$, 
and $t_3=\ui$.
We have $m_i(t_1)=-a(i)+a(i-1)$;
$$m_i(t_2)=\phi_i(\tf_i^{a(\ol{i+1})}b_i)+\sum_{t'\in S}\ipp{i}{\wt t'},
\hskip.6em S=\bigl\{\,\tf_{i+1}^{a(\ol{i+2})}b_{i+1},
\tf_{i+1}^{a(i+1)}b_{i+1},\tf_i^{a(i)}b_i,\tf_{i-1}^{a(i-1)}b_{i-1}\,\bigr\}$$ 
$$m_i(t_2)=-a(\ol{i+1})+a(\ol{i+2})+a(i+1)-2a(i)+a(i-1);$$
$m_i(t_3)=a(\ol{i})-a(\ol{i+1})+m_i(t_2)$.
Since $a(\ol{i})-a(\ol{i+1})\le0$, we have $m_i(t_3)\le m_i(t_2)$, so the
maximum of the $m_i(t)$ is $\max\bigl(m_i(t_1),m_i(t_2)\bigr)$, and (1) holds.
If $i=1$, the proof is similar, using $a(0)=0$. 

Suppose that $i=n-1$, $n>2$. Then 
$\phi_{n-1}(b)=\max\bigl(m_{n-1}(t_1),m_{n-1}(t_2),m_{n-1}(t_3)\bigr)$
where $t_1=\tf_{n-1}^{a(n-1)}b_{n-1}$, $t_2=\tf_{n-1}^{a(\ol{n})}b_{n-1}$,
$t_3=\ui$. Now $m_{n-1}(t_1)=-a(n-1)+a(n-2)$.
We have
$$m_{n-1}(t_2)=\phi_{n-1}(\tf_{n-1}^{a(\ol{n})}b_{n-1}) + 
\sum_{t'\in S}\ipp{i}{\wt t'},\quad 
S=\bigl\{\,\tf_{n}^{a(n)}b_{n},
\tf_{n-1}^{a(n-1)}b_{n-1},\tf_{n-2}^{a(n-2)}b_{n-2}\,\bigr\}.$$
Since
$\ipp{n-1}{\wt\tf_n^{a(n)}b_n}=-a(n)\ipp{n-1}{\alpha_n}=-a(n)(-2)=2a(n)$,
then
$$m_{n-1}(t_2)=-a(\ol n)+2a(n)-2a(n-1)+a(n-2).$$
We have $m_{n-1}(t_3)=a(\ol{n-1})-a(\ol n)+m_{n-1}(t_2)\le m_{n-1}(t_2)$, since
$a(\ol{n-1})\le a(\ol n)$. Thus $\phi_{n-1}(b)=
\max\bigl(m_{n-1}(t_1),m_{n-1}(t_2)\bigr)$, 
which proves (2). If $n=2$, the proof is similar.
The proof of (3) is straightforward.

\bigskip

\proclaim Lemma 4.3. Suppose that $b$ is the element as in Lemma 4.2. If
$1\le i<n-1$, then
$\tf_ib$ is gotten from $b$ by increasing $a(i)$ by 1 if
$a(\ol{i+2})-a(\ol{i+1})>a(i)-a(i+1)$ and increasing $a(\ol{i+1})$
by 1 otherwise; 
$\tf_{n-1}b$ is gotten from $b$ by increasing $a(n-1)$ by 1 if
$a(n-1)-a(n)\ge a(n)-a(\ol{n})$, and increasing $a(\ol{n})$
by 1 otherwise; 
$\tf_nb$ is gotten from $b$ by increasing $a(n)$ by 1.

\noindent Proof: Write $b=c\otimes d$ where 
$c=u_\infty\otimes \tf_1^{a(\ol2)}b_1
\otimes\cdots \tf_{n-1}^{a(\ol n)}b_{n-1}$,
$d=\tf_n^{a(n)}b_n\otimes\tf_{n-1}^{a({n-1})}b_{n-1}
\otimesc\tf_1^{a({1})}b_1$.
Compare $\phi_i(c)$ to $\epsilon_i(d)$.
Suppose that $1<i<n-1$. Then as in the proof of Lemma 4.2, we have
$\phi_i(c)=a(\ol{i+2})-a(\ol{i+1})$. Similarly, using (2.3),
$\epsilon_i(d)=a(i)-a(i+1)$.
Thus, if $\phi_i(c)>\ei(d)$, that is if $ a(\ol{i+2})-a(\ol{i+1})
>a(i)-a(i+1)$, then from (2.1), $\tf_ib=\tf_ic\otimes d$; $\tf_ic$ is gotten
from $c$ by increasing $a(\ol{i+1})$ by 1. If $\phi_i(c)\le\ei(d)$,
then $\tf_ib=c\otimes\tf_i d$, and $\tf_id$ is gotten from $d$ by
increasing $a(i)$ by 1. So the Lemma holds in this case.
We get the same result if $i=1<n-1$.
If $i=n-1$, then $\phi_{n-1}(c)=-a(\ol n)$, while
$\epsilon_{n-1}(d)=a(n-1)-2a(n)$. Adding $a(n)$, we compare
$a(n)-a(\ol n)$ to $a(n-1)-a(n)$. For $i=n$,
$\phi_n(c)=a(\ol n)$;
$\epsilon_n(d)=a(n)$.
Since $a(\ol n)\le a(n)$, then $\phi_n(c)\le\epsilon_n(d)$,
so $\tf_n(b)=c\otimes \tf_nd$, and we increase $a(n)$ by 1.

\proclaim Lemma 4.4. Let $T$ be a one-rowed semi-standard tableau.
If $1\le i<n$, then $\tf_iT$ is gotten from $T$ as follows:
if the number of $\ol{i+1}$'s is greater than
the number of $i+1$'s, change
the rightmost $\ol{i+1}$ to $\ol\imath$ or get 0 if there are no $\ol{i+1}$'s;
otherwise change the rightmost $i$ to $i+1$ or get 0 if there are no $i$'s.
If $i=n$, change the rightmost $n$ to $\ol n$ or get 0 if there are no
$n$'s.

\noindent Proof: Suppose that $1\le i <n$. 
Separate $T$ into 2 parts: the left part $T_l$ consists
of all numbers in $T$ which are $\le i+1$, and right part $T_r$
consists of all the other entries in $T$. Then $T=T_r\otimes T_l$.
We have $\phi_i(T_r)=$ the number of $\ol{i+1}$'s and
$\epsilon_i(T_l)=$ the number of $i+1$'s. Now use  (2.1).
The result is clear if $i=n$.

\proclaim Lemma 4.5. Let $T$ be a one-rowed semi-standard tableau.
If $1\le i<n$, let $r$ be the number of\/ $\ol{i+1}$'s in $T$,
let $s$ be the number of $i+1$'s, and let $t$ be the number of $i$'s.
Then $\phi_i(T)=\max(t,r-s+t)$. We have $\phi_n(T)=$ the number of $n$'s.

\noindent Proof: Suppose that $1\le i<n$. 
Compute $\phi_i(T)$ as $\max\{k:\tf_i^kT\ne0\}$.
Suppose that $r>s$; from
the previous Lemma, $\tf_iT$ is gotten by changing
an $i$ to $i+1$. Then the number of $i+1$'s is one greater; so 
$\tf_i^{r-s}T$ has the same number of $i+1$'s as $\ol{i+1}$s. Then
apply $\tf_i$ to $\tf_i^{r-s}T$ $t$ times, changing $i$ to $i+1$ until
we run out of $i$'s. Hence $\phi_i(T)=r-s+t$. If $r\le s$, applying
$\tf_i$ to $T$ changes an $i$ to $i+1$, so $\phi_i(T)=t$. It is clear
that $\phi_n(T)$ is the number of $n$'s.

\proclaim Proposition C.
Let $T$ be a semi-standard tableau with one
row and at least one 1. 
Let $a(x)$ be the number of entries in $T$ which are $\succ x$.
Let $b=F(T)=
u_\infty\otimes \tf_1^{a(\ol2)}b_1\otimes \tf_2^{a(\ol3)}b_2
\otimes\cdots\otimes \tf_{n-1}^{a(\ol n)}b_{n-1}
\otimes\tf_n^{a(n)}b_n\otimes\tf_{n-1}^{a({n-1})}b_{n-1}
\otimesc\tf_1^{a({1})}b_1$ in the crystal
$\Bi\otimes B_1\otimes B_2\otimesc B_n\otimes B_{n-1}\otimesc B_1$,
$n\ge2$.
Then \hfill\break
(1) $\phi_i(F(T))=\phi_i(T)$ for $1<i\le n$;\hfill\break
(2) $\phi_1(F(T))=\phi_1(T)\,-$ the number of entries in $T$;\hfill\break
(3) $\ei(F(T))=\epsilon_i(T),\quad 1\le i\le n$;\hfill\break
(4) $\tf_iF(T)=F(\tf_iT)$ if $\tf_iT\ne0$.

\noindent Proof:
Suppose that $1<i<n-1$. From Lemma 4.2,
$$\phi_i(b)=\max\bigl(-a(\ol{i+1})+a(\ol{i+2})+a(i+1)-2a(i)+a(i-1),a(i-1)-a(i)\bigr).\eqno(4.2)$$
We have 
$$\eqalign{a(\ol{i+2})-a(\ol{i+1})&=\hbox{number of $\ol{i+1}$'s}=r\cr
a(i)-a(i+1)&=\hbox{number of $i+1$'s}=s\cr
a(i-1)-a(i)&=\hbox{number of $i$'s}=t,\cr}\eqno(4.3)$$
so $\phi_i(b)=\max(r-s+t,t)=\phi_i(T)$ by Lemma 4.5.
If $i=1$, the argument is similar,
except that the term $a(0)$
is missing in (4.2) and
$a(0)$ is the number of entries $\succ0$ in (4.3). Then
$\phi_i(b)=\phi_i(T)-$the number of entries of $T$.
If $i=n-1$, $n>2$, then from Lemma 4.2,
$$\displaylines{\phi_{n-1}(b)=
\max\bigl(-a(\ol n)+2a(n)-2a(n-1)+a(n-2),
a(n-2)-a(n-1)\bigr),\cr
a(n)-a(n-1)=\hbox{number of $\ol n$'s}=r\cr
a(\ol n)-a(n)=\hbox{number of $n$'s}=s\cr
a(n-2)-a(n-1)=\hbox{number of $n-1$'s}=t\cr}$$
and again $\phi_{n-1}(b)=\max(r-s+t,t)=\phi_{n-1}(T)$ by Lemma 4.5.
If $n=2$, the proof is similar. 
Finally, $\phi_n(b)=a(n-1)-a(n)=\hbox{number of $n$'s}=\phi_n(T)$.
This proves (1) and (2); the proof of (3) is similar.

To prove (4), suppose that $\tf_iT\ne0$.
If $1\le i<n-1$, then  $\tf_iF(T)$ is obtained from $F(T)$,
by Lemma 4.3, by increasing $a(\ol{i+1})$ by 1 if 
$a(\ol{i+2})-a(\ol{i+1})>a(i)-a(i+1)$ and by increasing $a(i)$ by
1 otherwise. From Lemma 4.4, $\tf_iT$ is gotten from
$T$ by changing an $\ol{i+1}$ to $\bari$ if 
the number of $\ol{i+1}$'s is greater than the number of $i+1$'s.
Now the number of $\ol{i+1}$'s is equal to $a(\ol{i+2})-a(\ol{i+1})$,
the number of $i+1$'s is equal to $a(i)-a(i+1)$. Thus
(4) holds in this case. The argument for $i=n-1$ is similar, and 
for $i=n$ it is easy.
\bigskip
\proclaim Lemma C. Suppose that $T$ is a semi-standard tableau whose
first row consists entirely of 1's; then $\eo(T)=0$.

\noindent Proof: If $\te_1T\ne0$, then $\te_1T$ is obtained from $T$
by changing a 2 to a 1 or changing a $\ol1$ to a $\ol2$. Since 
the first row of $T$ is all $1$'s, if there were a $\ol1$
in $T$, then 1 and $\ol1$ would occur in the same column, and then
condition (4.1) would fail. Changing a 2 to a 1 would leave
two 1's in the same column, so $\te_1T$ would not be semi-standard, 
violating [KN 4.3.2]. So $\te_1T=0$.

\proclaim Theorem C.
Let $B(\infty)$ be the crystal basis of $U_q^-(\ger g)$
where $\ger g$ has type $C_n$.
For each $i$, $1\le i\le n$, define
$$B(i)=B_i\otimes B_{i+1}\otimesc B_{n-1}\otimes B_n\otimes B_{n-1}
\otimesc B_i.$$
Let $\Psi$ be the unique morphism from $B(\infty)$ to
$\Bi\otimes B(1)\otimes B(2)\otimesc B(n)$
which maps $u_\infty$ to 
$u_\infty\otimes(b_1\otimes b_2\otimes
\cdots\otimes b_n\otimes b_{n-1} \otimesc b_1)\otimes
(b_2\otimes b_3\otimes
\cdots b_n\otimes b_{n-1} \otimesc b_2)\otimesc
 b_n$.   
For $b\in B(\infty)$, suppose that $\lambda$ is a dominant weight such that
the tableau $T$ associated to $\bar\pi_\lambda(b)$ is large.
Let $a(i,j)$ be the number of entries in the $i$-th row of $T$ which
are $\succ j$; for each $i$, define
$$\displaylines{\beta_i=\tf_i^{a(i,\ol{i+1})}b_i
\otimes\tf_{i+1}^{a(i,\ol{i+2})}b_{i+1}
\otimes\cdots\otimes \tf_{n-1}^{a(i,\ol n)}b_{n-1}
\otimes\tf_n^{a(i,n)}b_n\otimes\tf_{n-1}^{a(i,n-1)}b_{n-1}
\otimesc\tf_i^{a(i,i)}b_i,\cr 
F(T)=\ui\otimes \beta_1\otimes\beta_2\otimesc\beta_n.\cr}$$
Then $\Psi(b)=F(T)$.

\noindent
Proof: The proof is essentially the same as that of Theorem A. We 
use induction on $n$, starting at $n=1$.
For $n=1$, the argument is the same as in the first paragraph of
the proof of Theorem A, except that the one-rowed tableau $T$ has
entries $1$ and $\ol1$ instead of $1$ and $2$. 
The proof then proceeds as for Theorem A, using
Proposition C  and Lemma C instead of Proposition A and Lemma A.

\medskip
\proclaim Corollary C. The image of $\Psi$ consists of all elements
 of the form
$\ui\otimes \beta_1\otimes\beta_2\otimesc\beta_n$ 
where
$$\beta_i=\tf_i^{a(i,\ol{i+1})}b_i\otimes\tf_{i+1}^{a(i,\ol{i+2})}b_{i+1}
\otimesc \tf_{n-1}^{a(i,\ol n)}b_{n-1}
\otimes\tf_n^{a(i,n)}b_n\otimes\tf_{n-1}^{a(i,{n-1})}b_{n-1}
\otimesc\tf_i^{a(i,i)}b_i$$ and $\{a(i,j):1\le i\le n, i\le j\le n\}$,
$\{a(i,\barj):1\le i\le n, i+1\le\barj\le n\}$ are sets of integers
such that
$$0\le a(i,\ol{i+1})\le a(i,\ol{i+2})\le\cdots\le a(i,\ol n)\le a(i,n)\le
a(i,n-1)\le\cdots\le a(i,i), 1\le i\le n.$$

\noindent Proof: It follows from Theorem C that the image of $\Psi$ is
a subset of the set of elements of the form given in the Corollary.
On the other hand, given sets $\{a(i,j)\}$ and $\{a(i,\barj)\}$
of integers with the given conditions, build a large tableau $T$ with
nondecreasing entries in each row, with $a(i,x)$ entries in row $i$
which are $\succ x$.
Condition
(4.1) holds, since $i$ and $\bar i$ are never in the same column
since $T$ is large. Then $F(T)=\Psi(b)$ for $b=\bar\pi_\lambda(b)$
for suitable $\lambda$, and the result follows.

\bigskip

\noindent{5. \bf Type B}
\medskip\noindent
In this section we asssume that $\gg$ is of type $B_n$. Let $\alpha_1,
\ldots \alpha_{n-1}$ be the short roots, and $\alpha_n$ the long
root. We have $\langle h_i, \alpha_{i+1}
\rangle=-1$ for $1\le i\le n-2$, $\langle h_{n},\alpha_{n-1}\rangle=-2$,
$\langle h_i,\alpha_j\rangle=0$ if $|i-j|>1$.
The irreducible representation $V(\Lambda_1)$
of highest weight $\Lambda_1$ has dimension $2n+1$ over
$\Qq$, and the elements of its crystal basis are \sbox{$i$}$\,$ and 
\sbox{$\ol \imath$}$\,$, $1\le i\le n$, and \sbox{0}$\,$.
 The crystal graph of $B(\Lambda_1)$ is
$$\sbox{1}\tto1\sbox{2}\tto2\cdots\tto{n-1}\sbox{$n$}
\tto n\sbox0
\tto n\sbox{$\ol n$}\tto{n-1}\cdots\tto2
\sbox{\lower1.5pt\hbox{$\ol2$}}\tto1
\sbox{\lower1pt\hbox{$\ol1$}}\,.$$
We use the order $\prec$ where 
$$1\prec2\prec\cdots\prec n\prec0
\prec \ol n\prec\ol{n-1}\prec\cdots\prec\ol1.$$

Let $\omega_i=\Lambda_i$ for $1\le i\le n-1$, and $\omega_n=2\Lambda_n$.
The irreducible module
 $V(\omega_k)$ 
(called the anti-symmetric tensor representation
with highest weight $\omega_k$ in [KN, p.~320]) can be embedded
in $V(\Lambda_1)^{\otimes k}$. We identify the column
$$\vcenter{\by $i_1$\cr $i_2$\cr$\vdots$\cr
 $i_k$\cr\ey} \qquad\hbox{with~~~}
\sbox{$i_1$}\otimes\sbox{$i_2$}\otimes\cdots\otimes\tb{$i_k$}\in
B(\Lambda_1)^{\otimes k},$$
where each $i_j$ is $m$, $\ol m$, $1\le m\le n$, or 0.
Define $I_k^{(B)}$ to be the set of
all columns as above, where $1\preceq i_1\preceq\cdots\preceq i_k\preceq\ol1$
and any element other than 0 cannot occur more than once,
and if $i_j=p$ and $i_l=\ol p$, then 
$j+(k-l+1)\le p$. Then [KL, Prop.~5.3.1] the crystal basis $B(\omega_k)$
is equal to $I_k^{(B)}$.

In [KN, 5.7], a dominant weight $\lambda=\sum_{i=1}^n\lambda_i\Lambda_i$
is said to have type $(E)$ if $\lambda_n$ is even, say $2m$. 
Then $V(\lambda)$ can be embedded in 
$\VL^{\otimes \lambda_1}\otimes V(\Lambda_2)^{\otimes \lambda_2}
\otimesc V(\Lambda_{n-1})^{\otimes\lambda_{n-1}}\otimes
V(\omega_n)^{\otimes m}$, so its crystal
basis $B(\lambda)$ can be viewed as a subcrystal of
$B(\Lambda_1)^{\otimes\lambda_1}\otimes B(\Lambda_2)^{\otimes\lambda_2}
\otimesc B(\omega_n)^{\otimes m}$.
Then $b\in B(\lambda)$ is identified with a Young tableau, using
the right to left convention on columns, as in types
$A$ and $C$; the columns are in 
$I_k^{(B)}$.

Note that $V(\omega_n)$ is not $V(\Lambda_n)$, but $V(2\Lambda_n)$; 
$V(\Lambda_n)$ is the spin representation.
We do not need to use the spin representation, since for our
purposes, we can choose $\lambda_n$ to be even.

We will call a tableau $T$ 
{\it almost semi-standard\/} if its columns
are in $I_k^{(B)}$, the entries in the rows are non-decreasing
in the order $\prec$, and 0 occurs at most once in each row.
 Such a tableau is called large as in Definition 3.1.
Lemma 4.1 holds, with the same proof.

A semi-standard $B$-tableau of type $(E)$
is by definition [KN, p.~327] one
which is almost semi-standard, and for which two conditions
[KN, M.N.1, M.N.2] hold:
the first concerns successive columns where
$i$ and $\ol\imath$ occur in one of the columns, and this cannot happen
for a large tableau; the second concerns two successive columns with
$i$ in the left column and $\bari$ in the right column, in a lower row.
This also cannot happen in a large tableau. So a large semi-standard
tableau is the same thing as a large, almost semi-standard tableau. 

The basic result we need 
for type $B$ [KN, Theorem 5.7.1] is that if $\lambda$
has type $(E)$ then  $B(\lambda)$
is the set of semi-standard type $(E)$
$B$-tableaux of shape $V(\lambda)$.
Suppose that $b\in B(\infty)_\xi$, where $\xi\in Q_-$; then
$\xi=-\sum_{i=1}^nk_i\alpha_i$ where each $k_i$ is a non-negative
integer. Pick $\lambda_i>k_i$, with $\lambda_n$ even, say 
$\lambda_n=2m$.
Lemma 3.2 still holds, with the same proof, so $\bar\pi_\lambda(b)$
is a large semi-standard tableau.

\proclaim Lemma 5.1. Let $b$ be the element $u_\infty\otimes \tf_1^{a(\ol2)}b_1
\otimes \tf_2^{a(\ol3)}b_2
\otimes\cdots \otimes\tf_{n-1}^{a(\ol n)}b_{n-1}
\otimes\tf_n^{a(n)}b_n\otimes\tf_{n-1}^{a({n-1})}b_{n-1}
\otimesc\tf_1^{a({1})}b_1$ of the crystal
$\Bi\otimes B_1\otimes\cdots\otimes B_n\otimesc B_1$,
 where $0\le a(\ol2)\le\cdots\le a(\ol n)
\le a(n)/2\le a({n-1})\le\cdots\le a(1)$, $n\ge2$. Let $a(0)=0$.
Then
\hb(1) for $1\le i<n-1$, $\phi_i(b)$ is as in Lemma 4.2(1);
\hb(2) $\phi_{n-1}(b)=\max\bigl(-a(\ol n)+a(n)-2a(n-1)+a(n-2),
a(n-2)-a(n-1)\bigr)$;
\hb(3) $\phi_n(b)=-a(n)+2a(n-1)$.

\noindent Proof: If $1\le i<n-1$, this is the same as Lemma 4.2(1).
Suppose that $i=n-1$. The argument is the same as for the proof of
Lemma 4.2(2), except that 
$\ipp{n-1}{\wt\tf_n^{a(n)}b_n}=-a(n)\ipp{n-1}{\alpha_n}=-a(n)(-1)=a(n)$
(instead of $2a(n)$ in type C.) 

For (3), if $i=n$, 
$\phi_n(b)=\max\bigl(m_n(t_1),m_n(t_2)\bigr)$ where
$t_1=\tf_n^{a(n)}b_n$ and $t_2=\ui$. Now
$$\eqalign{m_n(t_1)&=
\phi_n(\tf_n^{a(n)}b_n)+\ipp n{\wt\tf_{n-1}^{a(n-1)}b_{n-1}}\cr
&=-a(n)-a(n-1)\ipp n{\alpha_{n-1}}=-a(n)+2a(n-1).\cr}$$
Similarly,
$$m_n(t_2)=2a(\ol n)-2a(n)+2a(n-1).$$
Since $2a(\ol n)\le a(n)$ by hypothesis, then
$m_n(t_2)\le m_n(t_1)$, so $\phi_n(b)=m_n(t_1)$, and (3) is proved.

\proclaim Lemma 5.2. 
Let $b$ be the element of Lemma 5.1. For $1\le i<n-1$, $\tf_ib$
is gotten from $b$ as in Lemma 4.3; $\tf_{n-1}b$ is gotten from $b$
by increasing $a(n-1)$ by 1 if $a(n-1)-a(n)\ge-a(\ol n)$ and increasing
$a(\ol n)$ by 1 otherwise; $\tf_nb$ is gotten from $b$ by increasing
$a(n)$ by 1.

\noindent Proof: Write
$b=c\otimes d$ where 
$$c=u_\infty\otimes \tf_1^{a(\ol2)}b_1
\otimes\cdots\otimes \tf_{n-1}^{a(\ol n)}b_{n-1},\qquad
d=\tf_n^{a(n)}b_n\otimes\tf_{n-1}^{a({n-1})}b_{n-1}
\otimesc\tf_1^{a({1})}b_1.$$ 
If $1\le i<n-1$, the proof is 
as for Lemma 4.3.
For $i=n-1$, $\phi_{n-1}(c)=-a(\ol n)$. Using (2.3),
$\epsilon_{n-1}(d)=\epsilon_{n-1}(\tf_{n-1}^{a(n-1)}b_{n-1})
-\ipp{n-1}{\wt \tf_n^{a(n)}b_n}=a(n-1)-a(n)\ipp{n-1}{\alpha_n}=a(n-1)-a(n)$.

For $i=n$, $\phi_n(c)=\ipp n{\tf_{n-1}^{a(\ol n)}b_{n-1}}=
-a(\ol n)\ipp n{\wt \alpha_{n-1}}=2a(\ol n)$;
$\epsilon_n(d)=a(n)$.
Since $a(\ol n)\le a(n)/2$, then $\phi_n(c)\le\epsilon_n(d)$,
so $\tf_n(b)=c\otimes \tf_nd$, and we increase $a(n)$ by 1.

\proclaim Lemma 5.3. Let $T$ be a semi-standard tableau with one row.
If $1\le i<n$, $\tf_iT$ is obtained from $T$ as in Lemma 4.4;
$\tf_nT$ is obtained from $T$ by changing 0 to $\ol n$ if there is
a 0 in $T$, by changing the right-most $n$ to 0 if there is an $n$
but no 0
in $T$, and getting 0 otherwise.

\noindent Proof: This is similar to Lemma 4.4.

\proclaim Lemma 5.4. Let $T$ be a one-rowed semi-standard tableau.
If $1\le i<n$, let $r$ be the number of\/ $\ol{i+1}$'s in $T$,
let $s$ be the number of $i+1$'s, and let $t$ be the number of $i$'s.
Then $\phi_i(T)=\max(t,r-s+t)$; $\phi_n(T)=$ twice the number of $n$'s
plus the number (0 or 1) of 0's.

\noindent Proof: For $1\le i<n-1$, this is the same as Lemma 4.5; for
$i=n$ it is straightforward.

\proclaim Proposition B. Let $T$ be a semi-standard tableau with one
row and at least one 1.  For $x\ne n$,
let $a(x)$ be the number of entries in $T$ which are $\succ x$,
and let $a(n)$ be twice the number of 
entries $\succ 0$ plus the number (0 or 1) of $0$'s.
Let $b=F(T)=u_\infty\otimes \tf_1^{a(\ol2)}b_1\otimes \tf_2^{a(\ol3)}b_2
\otimes\cdots \otimes\tf_{n-1}^{a(\ol n)}b_{n-1}
\otimes\tf_n^{a(n)}b_n\otimes\tf_{n-1}^{a({n-1})}b_{n-1}
\otimesc\tf_1^{a({1})}b_1$ be in
$\Bi\otimes B_1\otimes B_2\otimes\cdots\otimes B_n\otimes B_{n-1}
\otimesc B_1$, $n\ge2$.
Then
\hb (1) $\phi_i(F(T))=\phi_i(T)$ for $1<i\le n$;
\hb (2) $\phi_1(F(T))=\phi_1(T)\,-$ the number of entries in $T$;
\hb (3) $\ei(F(T))=\epsilon_i(T),\quad 1\le i\le n$;
\hb (4) $\tf_iF(T)=F(\tf_iT)$ if $\tf_iT\ne0$. 

\noindent Proof: For $1\le i<n-1$, this is the same as 
Proposition C.

For $i=n-1$, $n>2$, then from Lemma 5.1(2), 
$$\phi_{n-1}(b)=
\max\bigl(-a(\ol n)+a(n)-2a(n-1)+a(n-2),
a(n-2)-a(n-1)\bigr).$$
Write $\#i$ for the number of $i$'s, and $\#\succ i$ for the number
of symbols $\succ i$. Define $r=\#\ol n, s=\#n, t=\#(n-1)$.
Since $a(n-1)-a(n-2)=\#(n-1)=t$ then the first of the two terms being maximized is
$$\eqalign{a(n)-a(\ol n)-2a(n-1)+a(n-2)&=a(n)-a(\ol n)-a(n-1)+\bigl(a(n-2)-a(n-1)\bigr)\cr
&=a(n)-a(\ol n)-a(n-1)+t.\cr}$$
$$\eqalign{a(n)-a(\ol n)-a(n-1)&=2(\#\succ n)+\#0-\#\succ\ol n-\#\succ( n-1)\cr
&=(\#\succ n-\#\succ\ol n)-(\#\succ( n-1)-\#\succ n-\#0)
=r-s.\cr}$$
Thus $\phi_{n-1}(b)=\max(t,r-s+t)=\phi_{n-1}(T)$ by Lemma 5.3. 
The argument is similar if $n=2$.

For $i=n$, $\phi_n(b)=2a(n-1)-a(n)=2\#\succ( n-1)-2\#\succ 0-\#0=2\#n+\#0$.
So (1) and (2) hold; (3) is similar. 

For (4), if $1\le i<n-1$, this is similar to Proposition C.
If $i=n-1$, $\tf_{n-1}F(T)$ is obtained from $F(T)$, by Lemma 5.2
by increasing $a(\ol n)$ by 1
if $-a(\ol n)>a(n-1)-a(n)$ and increasing $a(n-2)$ by 1 otherwise.
Since $a(n)=2\#\succ n+\#0$, then
$$-a(\ol n)>a(n-1)-a(n)\iff 
\#\succ n+\#0-\#\succ\ol n>\#\succ(n-1)-\#\succ n 
\iff \#\ol n>\#n.$$
If $\tf_{n-1}T\ne0$, then $\tf_{n-1}T$ is obtained from $T$,
as in Lemma 4.4, by changing an $\ol n$ to $\ol{n-1}$ if the number
of $\ol n$'s is greater than the number of $n$'s, and changing 
an $n-1$ to an $n$ otherwise. So (4) holds for $i=n-1$.
It follows from Lemmas 5.2 and 5.3 that (4) holds for $i=n$. 

\proclaim Theorem B.
Let $B(\infty)$ be the crystal basis of $U_q^-(\ger g)$
where $\ger g$ has type $B_n$.
For each $i$, $1\le i\le n$, define
$$B(i)=B_i\otimes B_{i+1}\otimesc B_{n-1}\otimes B_n\otimes B_{n-1}
\otimesc B_i.$$
Let $\Psi$ be the unique morphism from $B(\infty)$ to
$\Bi\otimes B(1)\otimes B(2)\otimesc B(n)$
which maps $u_\infty$ to 
$u_\infty\otimes(b_1\otimes b_2\otimes
\cdots\otimes b_n\otimes b_{n-1} \otimesc b_1)\otimes
(b_2\otimes b_3\otimes
\cdots b_n\otimes b_{n-1} \otimesc b_2)\otimesc
 b_n$. 
For $b\in B(\infty)$, suppose that $\lambda$ is a dominant weight
of type $(E)$ such that
the tableau $T$ associated to $\bar\pi_\lambda(b)$ is large.
Let $a(i,j)$ be the number of entries in the $i$-th row of $T$ which
are $\succ j$, if $j\ne n$, while $a(i,n)$ is defined to be twice
the number of entries in row $i$ which are $\succ n$ plus the number
(0 or 1) of 0's in row $i$.
For each $i$, define
$$\displaylines{\beta_i= \tf_i^{a(i,\ol{i+1})}b_i
\otimes\cdots\otimes \tf_{n-1}^{a(i,\ol n)}b_{n-1}
\otimes\tf_n^{a(i,n)}b_n\otimes\tf_{n-1}^{a(i,n-1)}b_{n-1}
\otimesc\tf_i^{a(i,i)}b_i,\cr 
F(T)=\ui\otimes \beta_1\otimes\beta_2\otimesc\beta_n.\cr}$$
Then $\Psi(b)=F(T)$.

\noindent
Proof: The proof is essentially the same as that of Theorem C, using
Proposition B instead of Proposition C; Lemma C still holds for type
$B$ tableaux. We start the induction at $n=1$.
(Type $B_1$ is the same as type $A_1$; but the crystal
graph of $V(\Lambda_1)$ for type $B_1$ is the same as 
the crystal graph of $V(2\Lambda_1)$ for
type $A_1$.)
 When $n=1$, $b=\tf^k(\ui)$;
$T$ is a one-rowed tableau whose entries are $1$'s, $\ol1$'s, and possibly
one $0$. Then $k$ is twice the number of $\ol1$'s plus the number of
$0$'s; since this is the same as $a(1,1)$, then $\Psi(b)=F(T)$.
(When $n=2$, $B_2$ is the same as $C_2$, and $V(\Lambda_1)$ for
type $B_2$ is the same as $V(\Lambda_2)$ for type $C_2$. 
This does not matter for our proof, since Proposition B still
holds for $n=2$.)  
\medskip
\proclaim Corollary B. The image of $\Psi$ consists of all elements
 of the form
$\ui\otimes \beta_1\otimes\beta_{2}\otimesc\beta_n$ 
where
$$\beta_i=\tf_i^{a(i,\ol{i+1})}b_i
\otimes\tf_{i+1}^{a(i,\ol{i+2})}b_{i+1}
\otimesc \tf_{n-1}^{a(i,\ol n)}b_{n-1}
\otimes\tf_n^{a(i,n)}b_n\otimes\tf_{n-1}^{a(i,{n-1})}b_{n-1}
\otimesc\tf_i^{a(i,i)}b_i$$ and $\{a(i,j):1\le i\le n, i\le j\le n\}$,
$\{a(i,\barj):1\le i\le n, i+1\le\barj\le n\}$ are sets of integers
such that
$$\displaylines{
0\le a(i,\ol{i+1})\le a(i,\ol{i+2})\le\cdots\le a(i,\ol n)\le a(i,n)/2\le
a(i,n-1)\le\cdots\le a(i,i),\cr
{}\qquad\qquad\qquad 1\le i\le n.\cr}$$

\noindent Proof: This is similar to type $C$.

\bigskip

\noindent{6. \bf Type D}
In this section we asssume that $\gg$ is of type $D_n$.
The irreducible representation $V(\Lambda_1)$
of highest weight $\Lambda_1$ has dimension $2n$ over
$\Qq$, and the elements of its crystal basis are \sbox{$i$}$\,$ and 
\sbox{$\ol \imath$}$\,$, $1\le i\le n$.  
The crystal graph of $B(\Lambda_1)$ is
$$\sbox{1}\tto1\sbox{2}\tto2\cdots\tto{n-2}
\sbbox{$n-1$}
\kern-.5em\matrix{&\smbox n&\cr
\raise.7ex\hbox{\kern.2em$\scr n-1$}\nearrow&&\searrow
\raise.7ex\hbox{\kern-.2em$\scr n$}\cr
\cr
\cr
{\scr n}\searrow&&\nearrow{\scr n-1}\cr
&\smbox{\ol n}&\cr}
\kern-.7em\sbbox{$\ol{n-1}$}
\tto{n-2}\sbbox{$\ol{n-2}$}\tto{n-3}\cdots\tto2
\sbox{\lower1.5pt\hbox{$\ol2$}}\tto1
\sbox{\lower1pt\hbox{$\ol1$}}\,.$$

Let $\omega_k=\Lambda_k$, $1\le k\le n-2$, and let $\omega_{n-1}=
\Lambda_{n-1}+\Lambda_n$.
Then $V(\omega_k)$ can be embedded in 
$V(\Lambda_1)^{\otimes k}$, $1\le k\le n-1$.
We identify columns with tensors as with
types $A$, $B$, $C$. We use the partial order
$$1\prec2\prec\cdots\prec n-1\prec{n\atop\ol n}
\prec\ol{n-1}\prec\cdots\prec\ol2\prec\ol1.$$
This partial order is obtained from the crystal graph of $V(\Lambda_1)$.
Note that $n$ and $\ol n$ are not related.

Let $I_k^{(D)}$ denote the set of columns of length $k$ ($1\le k<n$)
with entries of the form $i$ or $\bari$, $1\le i\le n$,
where  for each $j$ ($1\le j<k$)
we have $i_j\prec i_{j+1}$ or $(i_j,i_{j+1})=(n,\ol n)$
or $(\ol n,n)$, and if $i_j=p$ and $i_l=\ol p$, then 
$j+(k-l+1)\le p$. Then the crystal basis $B(\omega_k)$ is equal to
$I_k^{(D)}$ [KN, Prop.~6.3.2].
 
\proclaim Definition 6.1. A dominant weight is said to have type $(W0)$
if it is of the form
$$\sum_{i=1}^{n-1}m_i\omega_i=\sum_{i=1}^{n-2}m_i\Lambda_i+m_{n-1}
(\Lambda_{n-1}+\Lambda_n).$$

In [KN, section 6] dominant weights of type $D$ are divided into 4
classes, and what we have called type $(W0)$ are
a subclass of $(W1)$ of [KN]. 
For $\lambda$ of type $(W0)$, $V(\lambda)$ can
be embedded in 
$\VL^{\otimes m_1}\otimes V(\Lambda_2)^{\otimes m_2}
\otimesc V(\Lambda_{n-2})^{\otimes m_{n-2}}\otimes 
V(\omega_{n-1})^{\otimes m_{n-1}}$, so its crystal
basis $B(\lambda)$ can be viewed as a subcrystal of
$B(\Lambda_1)^{\otimes m_1}\otimes B(\Lambda_2)^{\otimes m_2}
\otimesc B(\Lambda_{n-2})^{\otimes m_{n-2}}\otimes 
B(\omega_{n-1})^{\otimes m_{n-1}}$. 
Elements of $B(\lambda)$ are identified with Young tableaux,
again using the right-to-left convention on columns. We only need columns
of length at most $n-1$, and we do not need spin representations.

We call a Young tableau almost semi-standard if its columns are in $I_k^{(D)}$
($1\le k<n)$ and the elements in its rows satisfy $i_j\preceq i_{j+1}$.
Note that this means that $n$ and $\ol n$ cannot both be in the same row.
We again define a Young tableau to be large as in Definition 3.1.
A Young tableau is called semi-standard in [KN, section 6] if it is
almost semi-standard, and if certain conditions hold; 
these conditions are vacuous for a large semi-standard tableau.

Given an element $b\in B(\infty)_\xi$ where $\xi=\sum-k_i\alpha_i$,
pick integers $m_i$, $1\le i\le n-1$, such that $m_i>k_i$, $1\le i\le n-1$,
and $m_{n-1}>k_n$. Let $\lambda=\sum_{i=1}^{n-2}m_i\Lambda_i+m_{n-1}
(\Lambda_{n-1}+\Lambda_n)$. Then $\lambda$ is of type $(W0)$ and
$\bar\pi_\lambda(b)$ is large, by the analogue of Lemma 3.2.
According to [KN, Theorem 6.7.1], the crystal basis $B(\lambda)$
of a type $(W0)$ weight $\lambda$ corresponds to all semi-standard
tableaux of shape $\lambda$.

\proclaim Lemma 6.2. Let $b$ be the element $u_\infty\otimes \tf_1^{a(\ol2)}b_1
\otimes \tf_2^{a(\ol3)}b_2
\otimes\cdots \tf_{n-2}^{a(\ol{n-1})}b_{n-2}
\otimes\tf_n^{a(n)}b_n\otimes\tf_{n-1}^{a({n-1})}b_{n-1}
\otimesc\tf_1^{a({1})}b_1$ of the crystal
$\Bi\otimes B_1\otimes B_2\otimes\cdots\otimes B_{n-2}
\otimes B_n\otimes B_{n-1}\otimesc B_1$,
 where $0\le a(\ol2)\le a(\ol3)\le\cdots\le a(\ol{n-1})\le
\min(a(n),a(n-1))$, $\max(a(n),a(n-1))\le a(n-2)\le a(n-3)\le\cdots
\le a(1)$, $n\ge3$. Let $a(0)=0$.
Then
\hb(1) for $1\le i<n-1$, $\phi_i(b)$ is as in Lemma 4.2(1);
\hb(2) $\phi_{n-2}(b)=
\max\bigl(a(n-3)-a(n-2),a(n)+a(n-1)-a(\ol{n-1})-2a(n-2)+a(n-3)\bigr)$;
\hb(3) $\phi_{n-1}(b)=a(n-2)-a(n-1)$;
\hb(4) $\phi_n(b)=a(n-2)-a(n)$.

\noindent Proof: For $1\le i\le n-3$, this is the same as Lemma 4.2.
For $i=n-2$, $n>3$,
$\phi_{n-2}(b)=\max\bigl(m_{n-2}(t_1),m_{n-2}(t_2),m_{n-2}(t_3)\bigr)$
where $t_1=\tf_{n-2}^{a(n-2)}b_{n-2}$, $t_2=\tf_{n-2}^{a(\ol{n-1})}b_{n-2}$,
$t_3=\ui$. Then $m_{n-2}(t_1)=-a(n-2)-a(n-3)$; $m_{n-2}(t_2)$ equals 
$$\eqalign{
-&a(\ol{n-1})+ \hbox{the sum of $\ipp{n-2}{\wt t'}$ where $t'=\tf_n^{a(n)}b_n,\tf_{n-1}^{a(n-1)}b_{n-1}$, and $\tf_{n-3}^{a(n-3)}b_{n-3}$}\cr
&=-a(\ol{n-1})+a(n)+a(n-1)-2a(n-2)+a(n-3).\cr}$$
We have $m_{n-2}(t_3)\le m_{n-2}(t_2)$, as in the proof of Lemma 4.2.
It follows that $\phi_{n-2}(b)=\max\bigl(m_{n-2}(t_1),m_{n-2}(t_2)\bigr)$,
and (1) holds in this case. The argument is similar for $n=3$.
For $i=n-1$ and $i=n$, the proof is as in Lemma 3.3.

\proclaim Lemma 6.3.
Let $b$ be the element of Lemma 6.2. For $1\le i<n-2$, then $\tf_ib$
is gotten from $b$ as in Lemma 4.3; 
$\tf_{n-2}b$ is gotten from $b$ by increasing $a(n-2)$ by 1 if
$a(n-2)-a(n-1)-a(n)+a(\ol{n-1})\ge0$ and increasing $a(\ol{n-1})$ by 1
otherwise;
$\tf_{n-1}b$ is gotten from $b$
by increasing $a(n-1)$ by 1; $\tf_nb$ is gotten from $b$ by increasing
$a(n)$ by 1.

\noindent Proof: Write $b=c\otimes d$ where 
$$c=u_\infty\otimes \tf_1^{a(\ol2)}b_1
\otimes\cdots\otimes \tf_{n-2}^{a(\ol{n-1})}b_{n-2},\qquad
d=\tf_n^{a(n)}b_n\otimes\tf_{n-1}^{a({n-1})}b_{n-1}
\otimesc\tf_1^{a({1})}b_1.$$ 
If $1\le i<n-2$, the proof is essentially the same as
for Lemma 4.3.
For $i=n-2$, $\phi_{n-2}(c)=-a(\ol{n-1})$,
$\epsilon_{n-2}(d)=a(n-2)-a(n-1)-a(n)$.
If $i=n-1$, $\phi_{n-1}(d)=a(\ol{n-2})$,
$\epsilon_{n-1}(d)=a(n-1)$; since $a(\ol{n-1})\le a(n)$, then
$\tf_{n-1}(c\otimes d)=c\otimes\tf_{n-1}d$. A similar argument
works for $i=n$.

\proclaim Lemma 6.4. Let $T$ be a semi-standard tableau with one row.
If $1\le i<n-2$, $\tf_iT$ is obtained from $T$ as in Lemma 4.4;
$\tf_{n-1}T$ is obtained from $T$ by changing
the right-most $\ol n$ to $\ol{n-1}$ if there
are any $\ol n$'s, changing the right-most $n-1$ to $n$ if
there are no $\ol n$'s and there are some $n-1$'s, and getting
0 if there are no $n-1$'s or $\ol n$'s; $\tf_nT$ is obtained
from $b$ as for $\tf_{n-1}b$, changing the roles of $\ol n$ and
$n$.

\noindent Proof: This follows from the definitions.

\proclaim Lemma 6.5. Let $T$ be a one-rowed semi-standard tableau.
If $1\le i<n-1$, let $r$ be the number of\/ $\ol{i+1}$'s in $T$,
let $s$ be the number of $i+1$'s, and let $t$ be the number of $i$'s.
Then $\phi_i(T)=\max(t,r-s+t)$; 
$\phi_{n-1}(T)=\hbox{number of $\ol n$'s}+\hbox{number of $n-1$'s}$;
$\phi_{n}(T)=\hbox{number of $n$'s}+\hbox{number of $n-1$'s}$.

\noindent Proof: This is similar to Lemma 4.5.
\proclaim Proposition D. 
Let $T$ be a semi-standard tableau with one
row. Let $a(x)$ be the number of entries in $T$ which are $\succ x$
if $x$ is not equal to $n-1$ or $n$, while $a(n-1)$ is defined
to be the number of entries $\succeq n$ and $a(n)$ is the number 
$\succeq \ol n$. 
Let $b=F(T)=u_\infty\otimes \tf_1^{a(\ol2)}b_1\otimes \tf_2^{a(\ol3)}b_2
\otimes\cdots \otimes\tf_{n-2}^{a(\ol{n-1})}b_{n-2}
\otimes\tf_n^{a(n)}b_n\otimes\tf_{n-1}^{a({n-1})}b_{n-1}
\otimesc\tf_1^{a({1})}b_1$ be in
$\Bi\otimes B_1\otimes\cdots\otimes B_{n-2}\otimes B_n\otimes
B_{n-1}\otimes\cdots\otimes B_1$, $n\ge3$.
Then
\hb (1) $\phi_i(F(T))=\phi_i(T)$ for $1<i\le n$;
\hb (2) $\phi_1(F(T))=\phi_1(T)\,-$ the number of entries in $T$;
\hb (3) $\ei(F(T))=\epsilon_i(T),\quad 1\le i\le n$;
\hb (4) $\tf_iF(T)=F(\tf_iT)$ if $\tf_iT\ne0$.

\noindent Proof: For $1\le i<n-2$, this is the same as Lemma 4.6. For
$i=n-2$, $n>3$, then from Lemma 6.2,
$$\phi_{n-2}(b)=\max\bigl(a(n-3)-a(n-2),a(n)+a(n-1)-a(\ol{n-1})-2a(n-2)+a(n-3)\bigr).$$
Now $a(n-3)-a(n-2)=\#(n-2)=t$. If there are
no $\ol n$'s in $T$, then $a(n)=\#\succ \ol{n-1}$,
$a(\ol{n-1})=\#\succ \ol{n-2}$,
$a(n)-a(\ol{n-1})=\#\ol{n-1}=r$,
$a(n-2)-a(n-1)=\#(n-1)$.
Then $\phi_{n-2}(b)=\phi_{n-2}(T)$. If there are no $n$'s in $T$, then
$a(n-1)-a(\ol{n-1})=\#\ol{n-1}=r, a(n-2)-a(n)=\#(n-1)=s$,
and again $\phi_{n-2}(b)=\phi_{n-2}(b)$. If $i={n-1}$, 
$\phi_{n-1}(b)=a(n-2)-a(n-1)=\#\ol n+\#(n-1)=\phi_{n-1}(T)$.
The proofs for $n=3$ and $i=n$ for are similar. 

To prove (4), this is the same as for type $C$ if $i\le n-2$. For $i=n-2$,
$\tf_{n-2}F(T)$ is obtained from $F(T)$ by increasing $a(\ol{n-1})$ by 1
if $a(\ol{n-1})>a(n-2)-a(n-1)-a(n)$, and increasing $a(n-2)$ by
1 otherwise. As shown in the proof of (1),  
$a(\ol{n-1})-a(n-2)+a(n-1)+a(n)=\#\ol{n-1}-\#(n-1)$.
From Lemma 6.4, $\tf_{n-2}T$ is obtained from $T$
by changing an $\ol{n-1}$ to an $\ol{n-2}$ if $\#\ol{n-1}>\#(n-1)$
and changing an $n-2$
to an $n-1$ otherwise. So (4) holds for $i=n-2$. If $i=n-1$,
$\tf_{n-1}F(T)$ is obtained from $F(T)$, by Lemma 6.3, by
increasing $a(n-1)$ by 1. If $\tf_{n-1}T\ne0$,
and there are no $n$'s in $T$ then $\tf_{n-1}T$
is obtained from $T$, by Lemma 6.4,
 by changing an $\ol n$ to $\ol{n-1}$; 
then since $a(n-1)=\#\succeq n$,
$a(n-1)$
goes up by 1, and $a(n)$ is unchanged. If there are no $\ol n$'s 
in $T$, $\tf_{n-1}T$ is obtained from $T$ by changing an $n-1$ to
$n$, and this increases $a(n)$ by 1 and leaves $a(n-1)$ unchanged.
So (4) holds if $i=n-1$; the argument for $i=n$ is similar.

\proclaim Theorem D.
Let $B(\infty)$ be the crystal basis of $U_q^-(\ger g)$
where $\ger g$ has type $D_n$, $n\ge 2$.
For each $i$, $1\le i\le n-2$, define
$$B(i)=B_i\otimes B_{i+1}\otimesc B_{n-2}\otimes B_n\otimes
 B_{n-1}\otimesc B_i.$$
Define $B(n-1)=B_n\otimes B_{n-1}$.
Let $\Psi$ be the unique morphsim from $\Bi$ to
$$\Bi\otimes B(1)\otimes B(2)\otimesc B(n-1)$$ which maps $u_\infty$ to
$u_\infty\otimes(b_1\otimes b_2\otimes
\cdots\otimes b_{n-2}\otimes b_n\otimes b_{n-1} \otimesc b_1)\otimes
(b_2\otimes b_3\otimes
\cdots b_{n-2}\otimes b_n\otimes b_{n-1} \otimesc b_2)\otimesc
( b_n\otimes b_{n-1})$.
For $b$ in $B(\infty)$,
suppose that $\lambda$ is a dominant weight of type $(W0)$
such that the tableau $\bar\pi_\lambda(b)$ is large.
Let $a(i,j)$ be the number of entries in the $i$-th row of $T$ which
are $\succ j$ if $j$ is not $n-1$ or $n$; let 
$a(i,n-1)$ be the number of entries in the $i$-th row 
which are $\succeq n$, and let
$a(i,n)$ be the number of entries in the $i$-th row $\succeq\ol n$.
For $1\le i\le n-2$ define
$$\beta_i=\tf_i^{a(i,\ol{i+1})}b_i
\otimes\tf_{i+1}^{a(i,\ol{i+2})}b_{i+1}\otimes
\cdots\otimes\tf_{n-2}^{a(i,\ol{n-1})}b_{n-2}
\otimes\ttf n{i,n}\otimes\tf_{n-1}^{a(i,n-1)}
b_{n-1}\otimes\cdots\otimes
\tf_i^{a(i,i)}b_i$$
and define $\beta_{n-1}=\tf_n^{a(n-1,n)}b_n\otimes
\tf_{n-1}^{a(n-1,n-1)}b_{n-1}$. Let
$$F(T)=\ui\otimes\beta_1\otimes\beta_2\otimesc\beta_{n-1}.$$
Then $\Psi(b)=F(T)$.

\noindent Proof:  Use induction on $n$, beginning the induction
at $n=2$ when we get a one-rowed tableau $T$ since our large
tableaux for type $D_n$ have $n-1$ rows.
When $n=2$, the root system $D_2$ is the same as $A_1\times A_1$.
Then $b\in B(\infty)$ has the form $b=\tf_2^k\tf_1^l\ui=\tf_1^l\tf_2^k\ui$
for some non-negative integers $k$ and $l$. Apply
 $\Psi:B(\infty)
\to B(\infty)\otimes B_2\otimes B_1$, giving
$\Psi(b)=\ui\otimes\tf_2^kb_2\otimes\tf_1^lb_1$. The tableau $T$ contains
$1$'s, $2$'s, and either $\ol1$'s or $\ol2$'s. Suppose that
$k>l$; then $T$ has $l$ $\ol1$'s and $k-l$ $2$'s. Then $k$ is the
number of $\ol1$'s plus the number of $2$'s, which is equal to
$a(1,1)$ by definition, and $l=a(1,2)$, so $\Psi(b)=F(T)$.
The argument is similar if $l\ge k$.
  
 The proof proceeds
as for the other types, using Proposition D and Lemma C. 
(When $n=3$, $D_3$ is the same
as $A_3$, and $V(\Lambda_1)$ for $D_3$ is the same as $V(2\Lambda_1)$ for
$A_3$;
Proposition
$D$ is still valid for $n=3$.)

\medskip
\proclaim Corollary D. The image of $\Psi$ consists of all elements
 of the form
$\ui\otimes \beta_1\otimes\beta_2\otimesc\beta_{n-1}$ 
where for $1\le i\le n-2$,
$$\beta_i=\tf_i^{a(i,\ol{i+1})}b_i\otimes\tf_{i+1}^{a(i,\ol{i+2})}b_{i+1}
\otimesc \tf_{n-2}^{a(i,\ol{n-1})}b_{n-2}
\otimes\tf_n^{a(i,n)}b_n\otimes\tf_{n-1}^{a(i,{n-1})}b_{n-1}
\otimesc\tf_i^{a(i,i)}b_i,$$
$\beta_{n-1}=\tf_n^{a(n-1,n)}b_n\otimes
\tf_{n-1}^{a(n-1,n-1)}b_{n-1}$
 and $\{a(i,j):1\le i\le n-1, i\le j\le n\}$,
$\{a(i,\barj):1\le i\le n-1, i+1\le\barj\le n-1\}$ are sets of 
non-negative integers
such that
$$\displaylines{
a(i,\ol{i+1})\le a(i,\ol{i+2})\le\cdots\le a(i,\ol{n-1})
\le\min\bigl(a(i,n-1),a(i,n)\bigr),\cr
\max\bigl(a(i,n-1),a(i,n)\bigr)\le
a(i,n-2)\le\cdots\le a(i,i),\cr
{}\qquad\qquad\qquad 1\le i\le n-1.\cr}$$

\noindent Proof: As for type C.

\bigskip
\centerline{REFERENCES}
\medskip\frenchspacing
\baselineskip=12pt
\parindent=0pt \parskip=6pt

[K1] M. Kashiwara, On crystal bases of the $q$-analogue of universal
enveloping algebras, Duke Math.~J. {\bf63} (1991), 465--516.

[K2] M. Kashiwara, Crystal base and Littelmann's refined Demazure
character formula, Duke Math.~J. {\bf71} (1993), 839--858.

[K3] M. Kashiwara, On crystal bases, in Representations of Groups, 155--197,
Proc. Canadian Math. Soc. Annual Seminar (Banff, AB, 1994), 
B. N. Allison \& G. H. Cliff
(eds.), CMS Conf. Proc {\bf16}, Amer. Math. Soc. 1995. 

[KN] M. Kashiwara and T. Nakashima, Crystal graphs for representations
of the $q$-analogue of classical Lie algebras, J.~Algebra {\bf 165} (1994),
295--345.

\bye